\documentclass[12pt,a4paper]{article}
\usepackage[top=2cm, bottom=2.5cm, left=1.6cm, right=1.6cm]{geometry}
\usepackage{amsmath,amssymb,amsfonts,amsxtra,graphics,graphicx,amsthm,bbm,array,subfig}
\usepackage{mathrsfs}
\usepackage{hyperref}
\usepackage{CJK}
\usepackage{color}
\usepackage{wrapfig}

\usepackage[T1]{fontenc}
\usepackage[utf8]{inputenc}
\usepackage{authblk}
\usepackage{dcolumn}
\usepackage[vcentermath]{youngtab}

\newcommand{\bZ}{\ensuremath{\mathbb{Z}}}

\newcommand{\frakg}{\ensuremath{\mathfrak{g}}}

\newcommand{\fraksl}{\ensuremath{\mathfrak{sl}}}

\newcommand{\cW}{\mathcal{W}}

\newcommand{\Tr}{\mbox{Tr}}

\def\beaa{\begin{eqnarray*}}
\def\eeaa{\end{eqnarray*}}
\def\bee{\begin{equation*}}
\def\eee{\end{equation*}}
\def\bea{\begin{eqnarray}}
\def\eea{\end{eqnarray}}
\def\be{\begin{equation}}
\def\ee{\end{equation}}
\def\ba{\begin{align}}
\def\ea{\end{align}}

\newcommand{\bem}{\begin{pmatrix}}
\newcommand{\eem}{\end{pmatrix}}

\def\={\;  = \;}
\def\+{\, + \,}

\def\wh{\widehat}
\def\bar{\overline}

\def\rt2{\sqrt{2}}

\begin{document}
\Yboxdim4pt
\title{
Colored HOMFLY-PT polynomials \\ that distinguish mutant knots}

\author[1,2]{Satoshi Nawata}
\author[3]{P. Ramadevi}
\author[3]{Vivek Kumar Singh}
\affil[1]{Faculty of Physics, University of Warsaw,  ul. Pasteura 5,  02-093, Warsaw, Poland}
\affil[2]{Max-Planck-Institut f\"ur Mathematik, Vivatsgasse 7, D-53111 Bonn, Germany}
\affil[3]{Department of Physics, Indian Institute of Technology Bombay,  India, 400076}
\date{}

  \maketitle

\abstract{We illustrate from the viewpoint of braiding operations on WZNW conformal blocks how colored HOMFLY-PT polynomials with multiplicity structure can detect mutations.  As an example, we explicitly evaluate the $\yng(2,1)$-colored HOMFLY-PT polynomials that distinguish a famous mutant pair, Kinoshita-Terasaka and Conway knot. \\}




\vspace{.5cm}

\section{Introduction}

\allowdisplaybreaks

A knot is an embedding of a circle in three-dimensional space up to ambient isotopy. Two knots are regarded topologically equivalent if one knot can be continuously deformed into the other knot without cutting itself. Knots are such complicated objects that, even with current technique in algebraic topology, it is hard to decide if two knots are topologically equivalent. To address this problem, a number of topological invariants of knots have been introduced. Among them, Jones has discovered a significantly powerful knot invariant, called a Jones polynomial \cite{Jones:1985dw}. It turned out that the Jones polynomial can be naturally understood in the framework of quantum physics. Witten has proposed a formulation of the Jones polynomial as an expectation value of Wilson loop in Chern-Simons theory \cite{Witten:1988hf}. Inspired by Witten's work, numerous quantum knot invariants have been constructed. In particular, Witten's formulation reveals the relation of quantum knot invariants with the two-dimensional WZNW model. By means of braiding operations on $\widehat\frakg_k$ WZNW conformal blocks, one can construct quantum knot invariants colored by a representation $R$ of $\frakg$ \cite{Ramadevi:1992dh}. Jones, HOMFLY-PT and Kauffman polynomials are indeed associated to $\fraksl(2)$, $\fraksl(N)$ and $\mathfrak{so}/\mathfrak{sp}(N)$ quantum knot invariants with the fundamental representation, respectively.

\begin{figure}[ht]
\begin{center}
\includegraphics[scale=.9]{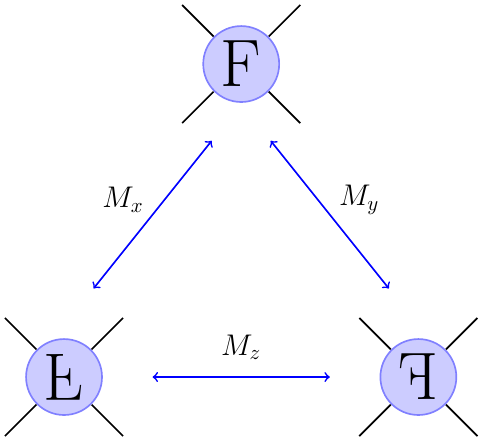}
\end{center}
\caption{Mutations on a two-tangle. Two out-going and two in-going orientations should  be added to four strands appropriately.}
\label{fig:basic-mutation}
\end{figure}

\begin{figure}[ht]
\begin{center}
\includegraphics[scale=.5]{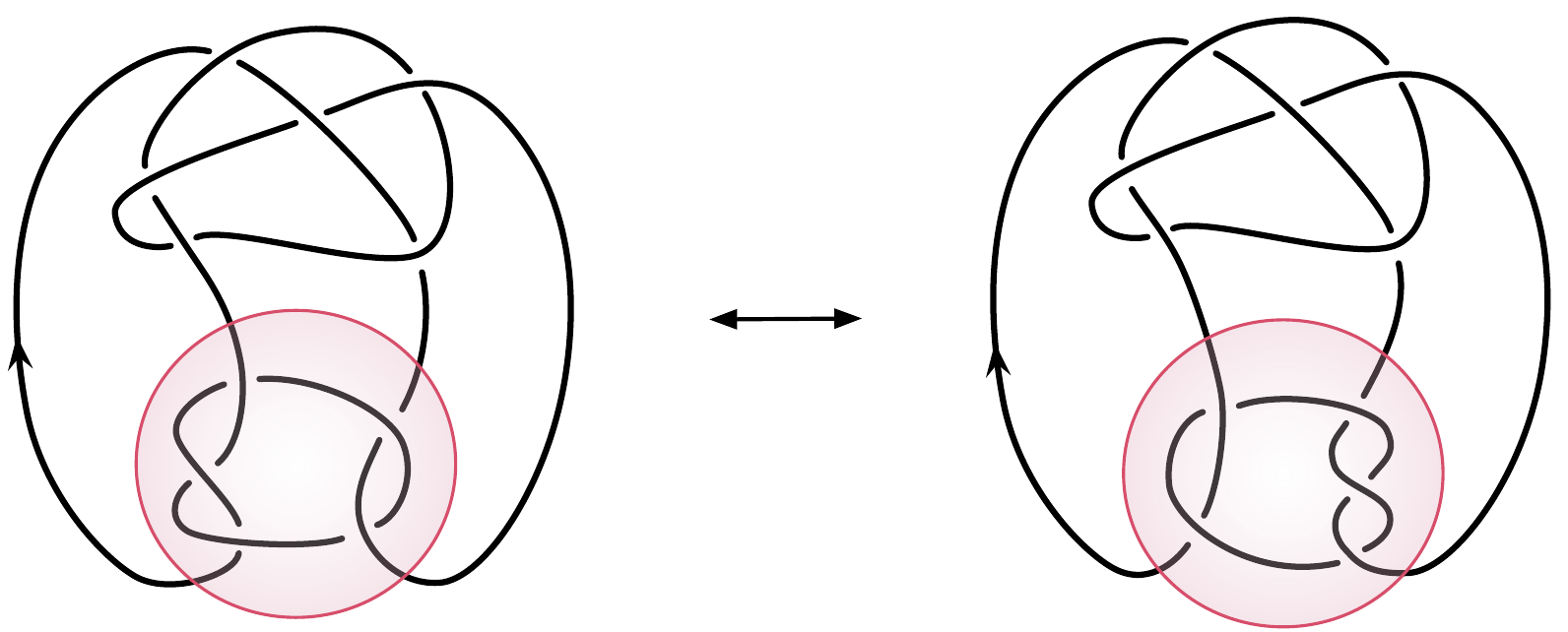}
\end{center}
\caption{The famous mutant pair: Kinoshita-Terasaka (left) and Conway (right) knots. The two knots are transformed into each other by the mutation $M_y$ (or $M_z$) on the shaded two-tangle region.}
\label{fig:KT-C}
\end{figure}

The quantum $\frakg$ knot invariants are very powerful and it is an important problem to understand whether they can distinguish mutant knots. A mutant pair \cite{Conway:1970} is obtained by performing a $180^\circ$ rotation about the horizontal axis ($M_x$), the vertical axis ($M_y$) or the axis perpendicular to the paper ($M_z$) on any two-tangle of a knot $K$. Note that the rotations are not independent so that $M_z=M_xM_y$ (Figure \ref{fig:basic-mutation}). It is well-known that many topological invariants cannot distinguish a mutant pair. In fact, a mutant pair shares the same uncolored Jones, HOMFLY-PT and Kauffman polynomials. Moreover, it is shown in \cite{Morton:1996} that quantum $\frakg$ invariants colored by a representation $R$ cannot distinguish mutant knots if a tensor product $R\otimes R$ is multiplicity-free, \textit{i.e.} the tensor product decomposes into irreducible representations with no repeated summands \cite[Theorem 5]{Morton:1996}. Among mutant pairs, Kinoshita-Terasaka $K_{KT}$ and Conway $K_C$ knots are the most famous example (Figure \ref{fig:KT-C}). Even though Khovanov homology can distinguish some mutant links, Kinoshita-Terasaka and Conway knots have the same Khovanov homology \cite{Wehrli:2003}. On the other hand, Morton and Cromwell have shown that  $\yng(2,1)$-colored HOMFLY-PT polynomials can detect this mutant pair by directly evaluating the difference of  invariants of  their satellites \cite{Morton:1996}. In fact, the tensor product of mixed (non-rectangular) representations such as $R=\yng(2,1)$ gives some irreducible representations more than once (multiplicity), and the multiplicity structure plays a pivotal role to detect mutant pairs. Therefore, taking \cite[Theorem 6]{Morton:1996} into account,  $\yng(2,1)$-colored quantum $\fraksl(N)$ invariants can distinguish the mutant pair only for $N>3$.

From the viewpoints of the cabling method \cite{Murakami:1989,Ochiai-Murakami} and the Reshetikhin-Turaev construction \cite{Morton:1998}, the reason why quantum knot invariants with multiplicity structure can detect a mutation is explained. On the other hand, the explanation has not been given yet directly in terms of WZNW conformal field theory although the approach from braiding operations on conformal blocks is equivalent to the Reshetikhin-Turaev construction \cite{Reshetikhin:1990} (Drinfel'd-Kohno theorem).\footnote{The explanation in \cite{Ramadevi:1994zb} does not deal with the multiplicity issue properly. Therefore, it is only applicable to multiplicity-free cases.} In this paper, we shall account for the reason why colored HOMFLY-PT polynomials with multiplicity structure can distinguish mutant knots from the viewpoint of WZNW model. Recently, Gu and Jockers have carried out impressive calculations of fusion matrices (quantum $6j$-symbols) for the representation $\yng(2,1)$ \cite{Gu:2014}. It turns out that the properties of quantum $6j$-symbols with multiplicity are different from multiplicity-free cases \cite{Kirillov:1989,Nawata:2013ppa}. In particular, signs called $3j$-phases play a crucial role for the detection. As an example, we explicitly compute the $\yng(2,1)$-colored HOMFLY-PT polynomials for Kinoshita-Terasaka and Conway knots by introducing three-boundary states with multiplicity indices in addition to the results in \cite{Gu:2014}.

The plan of the paper is as follows. In section 2, we will briefly discuss the relation between Chern-Simons theory, WZNW model and knot invariants with multiplicity. 
Then, realizing the mutations for $M_x$ and $M_y$ by braiding operations on a two-tangle, we will show how the multiplicity structure can detect the mutations. In section 3, we shall show that  the $\yng(2,1)$-colored HOMFLY-PT polynomials of  Kinoshita-Terasaka and Conway knots are indeed different by explicit computation.  We also attach a mathematica file for this computation to the arXiv page as an ancillary file.  
  
\section{Multiplicity can detect mutations}
To be self-contained, we shall briefly review the relation of Chern-Simons theory, WZNW model and quantum knot invariants. We refer the reader to \cite{Kohno:2002} for precise mathematical formulation and we mostly follow the notation of \cite{Gu:2014} in this paper. We consider $SU(N)$ Chern-Simons theory on $S^3$
\begin{equation}\label{CS-action}
S= {k \over 4\pi}\int_{S^3} \Tr \left(A \wedge dA + {2 \over 3} A\wedge A\wedge A\right)~
\end{equation}
where $k\in\bZ$ is a Chern-Simons level. Witten has formally constructed a quantum knot invariant as the expectation value of a Wilson loop $W_R(K)= \Tr_R [ P \exp\oint A]$  along a knot $K$ carrying representation $R$
\begin{equation}\label{Wilson}
\cW_R(K):=\langle W_R(K) \rangle= {\int {\cal D} \!A ~~e^{iS} W_R(K)\over \int {\cal D}\! A ~~ e^{iS}}~.
\end{equation}
Since the action \eqref{CS-action} is independent of the metric of $S^3$, it is a typical example of topological quantum field theories of Schwarz type. Therefore, although the expression \eqref{Wilson} is given by a Feynman integral over an infinite dimensional moduli space of gauge connections on $S^3$, the techniques of topological quantum field theories can be used to evaluate \eqref{Wilson} {exactly}. To this end, we decompose the three-sphere with a Wilson loop into a collection of three-manifolds with $S^2$ boundaries with marked points due to the Wilson loop. Then, the Chern-Simons partition function on a three-manifold with boundary $(S^2,p_1,\cdots,p_n)$ turns out to be an element of quantum Hilbert space on the boundary, which is isomorphic to the space of $n$-point $\widehat{\mathfrak{sl}}(N)_k$ WZNW conformal blocks \cite{Witten:1988hf}. Using this relation, the expectation value of a Wilson loop can be evaluated by braiding and fusion operations on $\widehat{\mathfrak{sl}}(N)_k$ WZNW conformal blocks. The outcome is indeed a Laurent polynomial with respect to $q=\exp\left( \frac{2\pi i}{k+N}\right)$ and the substitution $q^N=a$ leads to a two-variable knot invariant, a colored HOMFLY-PT polynomial of a knot.

For braiding operations, the basic building blocks are four-point conformal blocks. It is well-known that there are two bases for four-point conformal blocks that can be schematically depicted as follows:
\begin{eqnarray}\label{briad-eigen}
 {\raisebox{-1.0cm}{\includegraphics[width=3.5cm]{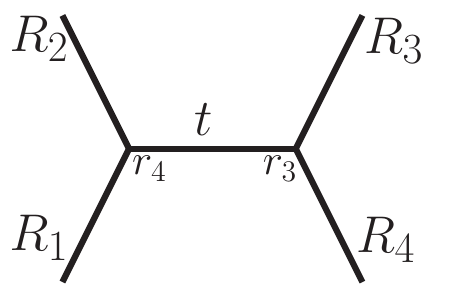}}}\!\!\!\!=|\phi^{(1)}_{t,r_3r_4}(R_1,R_2,R_3,R_4)\rangle \ ,&\quad&{\raisebox{-1.1cm}{\includegraphics[width=2.7cm]{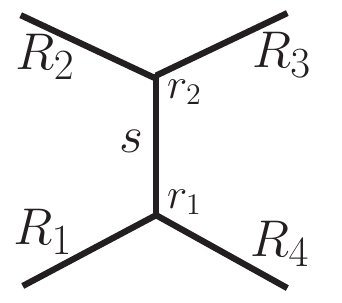}}}\!\!\!\!=|\phi^{(2)}_{s,r_1r_2}(R_1,R_2,R_3,R_4)\rangle\cr
 &&
\end{eqnarray}
Here, the WZNW primary fields are labelled by representations $R_i$ of $\wh\fraksl(N)_k$ and the intermediate states obey the fusion rule, \textit{i.e.} $t \in (R_1\otimes R_2)\cap(\bar{R}_3\otimes\bar{R}_4)$ in the first base $|\phi^{(1)}_{t,r_3r_4}(R_1,R_2,R_3,R_4)\rangle$ and $s\in (R_2\otimes R_3) \cap (\bar{R}_1\otimes\bar{R}_4)$ in the second base $|\phi^{(2)}_{s,r_1r_2}(R_1,R_2,R_3,R_4)\rangle$. Since  the decomposition for a tensor product of two  $\wh\fraksl(N)_k$ representations generally contains repeated summands of an irreducible representation, we use $r_i$ for its multiplicity label. The choice of basis  depends  on the braiding operator $b_i$ (half-monodromy).  In fact, the first basis $|\phi^{(1)}_{t,r_3,r_4}\rangle$ is an eigenstate of braiding operators $b_1^{(\pm)}$ ($b_3^{(\pm)}$) acting on the left (right) two strands whereas the second basis $|\phi^{(2)}_{s,r_1,r_2}\rangle$  is an eigenstate of a braiding operator $b_2^{(\pm)}$ acting on the middle two strands. Note that  the superscript $(+)$ is for parallel strands, or $(-)$ is for anti-parallel strands. The eigenvalue $\lambda^{(\pm)}_{R_i,R_j; R_kr_\ell}$ of a braiding operator $b_s^{(\pm)}$  depends on the two representations $R_i, R_j$ before the fusion, the intermediate representation $R_k$ and the multiplicity label $r_\ell$ after the fusion:
$$
\lambda^{(\pm)}_{R_i,R_j; R_k r_\ell} = \{R_i, R_j, \bar{R}_k, r_\ell\} q^{\pm(C_{R_i} + C_{R_j}-C_{R_k})/2}  
$$
where $C_R$ denotes the quadratic Casimir of a representation $R$.\footnote{Here we choose the canonical framing rather than the vertical framing used in \cite{Gu:2014}.} In addition, the phases $\{R_i, R_j, \bar{R}_k,r_l\}=\pm 1$ are called $3j$-phases that appear when  the two representations $R_i$ and $R_j$ are exchanged in  the Clebsch-Gordon coefficient:
$$
	\langle r_l R_k m_k| R_i m_i, R_j m_j \rangle = \{R_i, R_j, \bar{R}_k, r_l\} \langle r_l R_k m_k | R_j m_j , R_i m_i \rangle  \ .
$$
Furthermore, the first basis is transformed to the second one by the fusion/crossing matrices \cite{Kirillov:1989}
$$
	|\phi^{(1)}_{t,r_3r_4}(R_1, R_2, R_3, R_4) \rangle = \sum_{s,r_1,r_2} a^{t,r_3r_4}_{s,r_1r_2}\begin{bmatrix} R_1 & R_2 \\ R_3 & R_4 \end{bmatrix} | \phi^{(2)}_{s,r_1r_2}(R_1, R_2, R_3, R_4) \rangle \ .
$$
In fact, they satisfy the unitarity 
$$
	\sum_{t,r_3,r_4} a^{t,r_3r_4}_{s,r_1r_2}\begin{bmatrix} R_1 & R_2 \\ R_3 & R_4 \end{bmatrix} a^{t,r_3r_4}_{s',r'_1r'_2}\begin{bmatrix} R_1 & R_2 \\ R_3 & R_4 \end{bmatrix}^* = \delta_{s,s'} \delta_{r_1,r'_1} \delta_{r_2,r'_2} \ ,
$$
so that we also have the inverse relationship
$$
	|\phi^{(2)}_{s,r_1r_2}(R_1, R_2, R_3, R_4) \rangle = \sum_{t,r_3,r_4} a^{t,r_3r_4}_{s,r_1r_2}\begin{bmatrix} R_1 & R_2 \\ R_3 & R_4 \end{bmatrix}^* | \phi^{(1)}_{t,r_3r_4}(R_1, R_2, R_3, R_4) \rangle \ .
$$

With this setup, let us consider how the multiplicity structure can be used to distinguish mutant knots. 
In the following, we restrict ourselves to the case of either the representation $R=\yng(2,1)$ or its conjugate representation $\overline R=\overline{\yng(2,1)}$. It is straightforward to generalize the argument for more general settings. The tensor product of $R\otimes \bar{R}=(21;0)\otimes (0;21)$ and $R\otimes {R}=(21;0)\otimes (21;0)$ decomposes as follows:
\bea\label{irrep-[2,1]}
(21;0)\otimes (0;21) &=& (0;0)_0\oplus (1;1)_0\oplus (1;1)_1\oplus (2;2)_0 \oplus (2;1^2)_0  \\
&~&\oplus(1^2;2)_0 \oplus (1^2;1^2)_0 \oplus (21;21)_0  ,\nonumber \\
(21;0)\otimes (21;0)&=& 
(42;0)_0\oplus (2^3;0)_0 \oplus (31^3;0)_0 \oplus (321;0)_0\oplus(321;0)_1 \label {irrep-[2,1]-2}\\
&~&\oplus (41^2;0)_0 \oplus (3^2;0)_0\oplus (2^21^2;0)_0\nonumber~,
\eea
where the subscripts 0,1 are multiplicity indices and the irreducible representations $(1;1)$ and $(321;0)$ appear twice in $R\otimes \bar{R}$ and $R\otimes R$, respectively.\footnote{Following \cite{Gu:2014}, we use the composite labelling of the partition associated to an irreducible representation. The partition label $\lambda  = \lambda_1 \geqslant \lambda_2 \geqslant \cdots \geqslant \lambda_{N-1} \geqslant \lambda_N = 0$ can be recast in a composite manner,
\begin{align*}
	(\lambda) &= (\lambda_1, \lambda_2, \cdots, \lambda_N) \\
	&=( \mu_1, \mu_2 , \cdots, \mu_p, 0, \cdots, 0, -\nu_q, -\nu_{q-1}, \cdots, -\nu_1  ) \\
	&=( \mu; \nu )		\	,
\end{align*}
where $p+q \leqslant N$. Here the second line is obtained by subtracting the same integer from each $\lambda_i$. For example, the partition $ (4,3,2^{N-4},1,0)$ is equivalent to $(2,1,0^{N-4},-1,-2)=(21;21)$ by subtracting two from each row. Therefore, \eqref{irrep-[2,1]} can be expressed in the notation of partitions as follows:
\be
(2,1)\otimes (2^{N-2},1)=\varnothing\oplus (2,1^{N-2})\oplus (2,1^{N-2})\oplus (4,2^{N-2}) \oplus (3,1^{N-3})\oplus (4^2,2^{N-3}) \oplus (2^2,1^{N-4})\oplus (4,3,2^{N-4},1)\nonumber
\ee
} In these irreducible decompositions, one can take the $3j$-phases as in Table \ref{tab:3j} \cite{Gu:2014}. Importantly, the $3j$-phases for the two identical irreducible representations $(1;1)_0$ and $(1;1)_1$ ($ (321;0)_0$ and $ (321;0)_1$)  are different, which  plays a crucial role to detect a mutant pair.

\begin{table}[ht]\centering
\begin{tabular}{c|cccccccc}
$(t)_r$&$(0;0)_0$&$ (1;1)_0$&$ (1;1)_1$&$ (2;2)_0 $&$ (2;1^2)_0 $&$ (1^2;2)_0 $&$ (1^2;1^2)_0 $&$ (21;21)_0$\\ \hline
$\raisebox{\depth}      \{R,\bar{R},\bar{t},r  \}$&$+1$&$+1$&$-1$&$+1$&$-1$&$-1$&$+1$&$+1$\\
\end{tabular}
\end{table}
\begin{table}[ht]\centering

\begin{small}
\begin{tabular}{c|cccccccc}
$(t)_r$&$ (42;0)_0$&$ (2^3;0)_0 $&$ (31^3;0)_0 $&$ (321;0)_0$&$ (321;0)_1$&$ (41^2;0)_0 $&$ (3^2;0)_0$&$ (2^21^2;0)_0$\\ \hline
${\raisebox{\depth} \{R,{R},\bar{t},r\}}$&$+1$&$+1$&$+1$&$+1$&$-1$&$-1$&$-1$&$-1$\\
\end{tabular}
\end{small}
\caption{$3j$-phases for the tensor products of $R\otimes \bar{R}$ (upper) and $R\otimes R$ (lower)}  \label{tab:3j}
\end{table}

\begin{figure}[ht]
\center
\subfloat[Braiding operation for mutation $M_x$]{ \raisebox{2.3cm}{\includegraphics[width=0.36\textwidth]{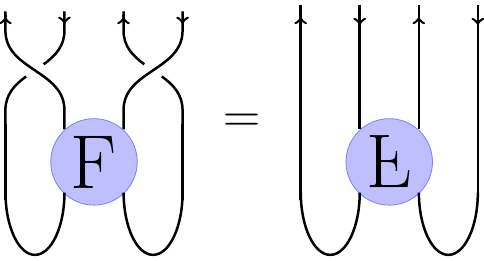}}} \quad \qquad  
\subfloat[Braiding operation for mutation  $M_y$] {\includegraphics[width=0.36\textwidth]{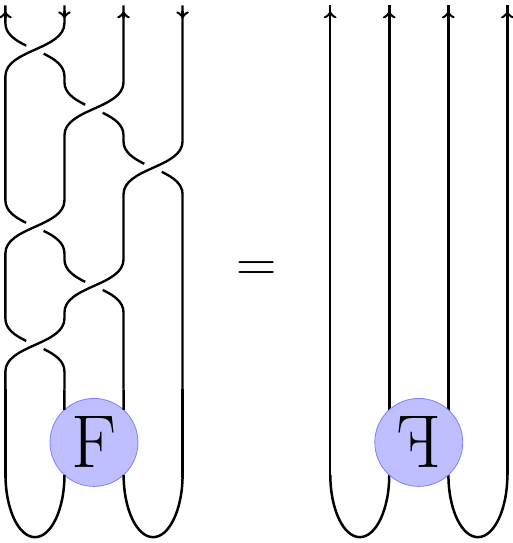}}
\caption{Braiding operations for mutations}\label{fig:braid-mutation}
\end{figure}

As shown in Figure \ref{fig:braid-mutation}(a), the mutation $M_x$ is realized by acting the braid word $b_1b_3^{-1}$ on a state $|\textrm{\textbf{F}}\rangle$ defined by the (4,0)-tangle. Using the braiding eigenvalues \eqref{briad-eigen}, it is easy to see that
\bea\label{Rx-braid}
|{ \raisebox{\depth}{\scalebox{1}[-1]{\textbf{F}}}}\rangle &=& b_1^{(-)} [b_3^{(-)}]^{-1} |\textrm{\textbf{F}}\rangle\\
 &=& \sum_{t,r_1,r_2}\{R,\bar{R},\bar{t},r_1\} \{R,\bar{R}, \bar t,r_2\} |\phi^{(1)}_{t,r_1,r_2}(R,\bar R,R,\bar R)\rangle\langle\phi^{(1)}_{t,r_1,r_2}(R,\bar R,R,\bar R)|\textrm{\textbf{F}}\rangle\nonumber~.
\eea
If $\langle\phi^{(1)}_{t,r_1,r_2}(R,\bar R,R,\bar R)|\textrm{\textbf{F}}\rangle\neq0$ for $(t,r_1,r_2)=((1;1),0,1)$ or $((1;1),1,0)$ , then we have $|\textrm{\textbf{F}}\rangle \neq |{ \raisebox{\depth}{\scalebox{1}[-1]{\textbf{F}}}}\rangle$ so that the $\yng(2,1)$-colored HOMFLY-PT polynomial can detect the mutation $M_x$. However, the special case $\langle\phi^{(1)}_{t,r_1,r_2}(R,\bar R,R,\bar R)|\textrm{\textbf{F}}\rangle=0$ for $(t,r_1,r_2)=((1;1),0,1)$ or $((1;1),1,0)$  occurs when the two-tangle $|\textrm{\textbf{F}}\rangle$ has a certain symmetry.
In fact, the direct computation using the three-boundary state \eqref{3-bdry} shows that $\langle\phi^{(1)}_{t,r_1,r_2}(R,\bar R,R,\bar R)|\textrm{\textbf{F}}\rangle$ vanishes for $(t,r_1,r_2)=((1;1),0,1)$ or $((1;1),1,0)$  in the example considered in \cite[Figure 2]{Morton:2009}. In addition, it is easy to see that  the two $3j$-symbols in \eqref{Rx-braid} cancel in the multiplicity-free case $r_1=r_2$ so that the braiding operation behaves as the identity operation.

In a similar manner, the braid word $[b_1^{(-)}]^{-1} b_2^{(+)}  [b_3^{(-)}]^{-1} [b_1^{(-)}]^{-1} b_2^{(+)} [b_1^{(-)}]^{-1}$ brings about the mutation $M_y$ (Figure \ref{fig:braid-mutation}(b)), which transforms  the state $|\textrm{\textbf{F}}\rangle$ to the state $|\reflectbox{\textbf{F}}\rangle$. Thus, we have
\bea\label{Ry-braid}
|\reflectbox{\textbf{F}}\rangle&=& \left([b_1^{(-)}]^{-1} b_2^{(+)}[b_1^{(-)}]^{-1} \right) b_1^{(-)} [b_3^{(-)}]^{-1}\left([b_1^{(-)}]^{-1} b_2^{(+)} [b_1^{(-)}]^{-1}\right)|\textrm{\textbf{F}}\rangle\cr
 &=& \sum_{t,s,r_1,r_2,r_3,r_4}  \{R,\bar{R},\bar{t},r_1\}\{R,\bar{R},\bar{t},r_2\}a_{s;r_3,r_4}^{u;r_6,r_5}\begin{bmatrix}
 R & \bar{R}\\
R& \bar{R}
\end{bmatrix}a_{t;r_1,r_2}^{s;r_4,r_3}\begin{bmatrix}
 R & \bar{R}\\
R& \bar{R}
\end{bmatrix}\times\cr
&&\hspace{2cm}|\phi^{(1)}_{t,r_1,r_2}(R,\bar R,R,\bar R)\rangle\langle\phi^{(1)}_{t,r_1,r_2}(R,\bar R,R,\bar R)|\textrm{\textbf{F}}\rangle\cr
 &=&\sum_{t,r_1,r_2}  \{R,\bar{R},\bar{t},r_1\}\{R,\bar{R},\bar{t},r_2\}|\phi^{(1)}_{t,r_2,r_1}(R,\bar R,R,\bar R)\rangle\langle\phi^{(1)}_{t,r_1,r_2}(R,\bar R,R,\bar R)|\textrm{\textbf{F}}\rangle~.
\eea
From the first to the second line, we apply the Racah backcoupling rule of the fusion matrices into the brackets of the first line
\bea\label{identity}
&&\{R,\bar{R},\bar{t},r_1\}\{R,\bar{R},\bar{t},r_2\}~a_{t;r_1,r_2}^{s;r_6,r_5}\begin{bmatrix}
 R & \bar{R}\\
R& \bar{R}
\end{bmatrix} \\
&&\hspace{2cm}=\sum_{u,r_3,r_4} (\lambda_{s;r_{6}}^{(-)})^{-1} ~a_{u;r_3,r_4}^{s;r_5,r_6}\begin{bmatrix}
 \bar{R}&R\\
R& \bar{R}
\end{bmatrix}^\ast ~  (\lambda_{u;r_{4}}^{(+)}) ~a_{u;r_3,r_4}^{t;r_1,r_2}\begin{bmatrix}
 \bar{R}&R\\
R& \bar{R}
\end{bmatrix} ~ (\lambda_{t;r_{2}}^{(-)})^{-1}\nonumber~,
\eea
and from the second to the third line, we have performed the summation over $s,r_3,r_4$, which results in the exchange of the multiplicity indices $r_1$ and $r_2$.
Again, we conclude that  $|\reflectbox{\textbf{F}} \rangle \neq |\textrm{\textbf{F}}\rangle$ unless $\langle\phi^{(1)}_{t,r_1,r_2}(R,\bar R,R,\bar R)|\textrm{\textbf{F}}\rangle=0$ for $(t,r_1,r_2)=((1;1),0,1)$ or $((1;1),1,0)$ in \eqref{Ry-braid}. In fact, Kinoshita-Terasaka and Conway knots are transformed by the mutation $M_y$ (Figure \ref{fig:KT-C} and Figure \ref{F-G}) and this is the reason why $\yng(2,1)$-colored HOMFLY-PT  polynomials distinguish them.

It is easy to see that the mutation $M_z$ generically does not behave as the identity operation in the multiplicity case since it can be realized as the composition of the previous two mutations $M_z=M_xM_y$. Besides, if the side two strands are with parallel orientation in Figure \ref{fig:braid-mutation}, we only need to consider the tensor product of $R\otimes R$ in \eqref{irrep-[2,1]-2}. Hence, a similar argument can be applied for the detection of mutations and we just avoid the repetition.


\subsection*{2-tangles}
As an illustrative example, we shall consider a mutation on a 2-tangle. Let us write a state $|\textrm{\textbf{F}}\rangle$ of a 2-tangle, and the state $ |{ \raisebox{\depth}{\scalebox{1}[-1]{\textbf{F}}}}\rangle$ and $ |\reflectbox{\textbf{F}}\rangle$ corresponding to the $M_x$--mutant and $M_y$--mutant tangle, respectively,   as
\begin{align}\label{F-state}
{\raisebox{-0.9cm}{\includegraphics[width=2cm]{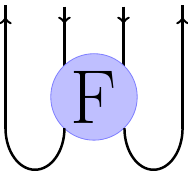}}} =\sum_{s,r_1,r_2} f_{s,r_1,r_2}|\phi^{(1)}_{s,r_1,r_2}(R,\bar{R},\bar{R},R)\rangle~,\\
{\raisebox{-0.9cm}{\includegraphics[width=2cm]{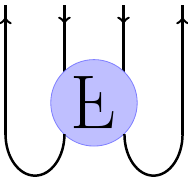}}} =\sum_{s,r_1,r_2}  {f}^{(x)}_{s,r_1,r_2}|\phi^{(1)}_{s,r_1,r_2}(R,\bar{R},\bar{R},R)\rangle ~,\\
{\raisebox{-0.9cm}{\includegraphics[width=2cm]{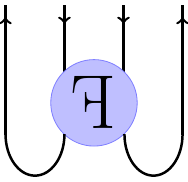}}} =\sum_{s,r_1,r_2}  {f}^{(y)}_{s,r_1,r_2}|\phi^{(1)}_{s,r_1,r_2}(R,\bar{R}.\bar{R},R)\rangle ~.
\end{align}
As we have seen, the coefficients are related by \eqref{Rx-braid} and \eqref{Ry-braid}, namely  
$$ f^{(x)}_{s,r_1,r_2}= (-1)^{r_1+r_2} f_{s,r_1,r_2}~,\qquad  f^{(y)}_{s,r_1,r_2}= (-1)^{r_1+r_2} f_{s,r_2,r_1}~.$$ 
Hence, it is easy to see that only representations $s$ with multiplicity are relevant, and in this situation $s=(1;1)$. The rest of the coefficients are equal. Hence the difference between tangle F and mutant of tangle $F$ will be
\begin{align}
|\textrm{\textbf{F}}\rangle- |{ \raisebox{\depth}{\scalebox{1}[-1]{\textbf{F}}}}\rangle&= 2f_{(1;1),1,0}|\phi^{(1)}_{(1;1),1,0}(R,\bar{R},\bar{R},R)\rangle+ 2f_{(1;1),0,1}|\phi^{(1)}_{(1;1),0,1}(R,\bar{R},\bar{R},R)\rangle~,\cr
|\textrm{\textbf{F}}\rangle- |\reflectbox{\textbf{F}}\rangle&= (f_{(1;1),0,1}+ f_{(1;1),1,0})\sum_{r_1 \neq r_2}|\phi^{(1)}_{(1;1),r_1,r_2}(R,\bar{R},\bar{R},R)\rangle~.\nonumber
\end{align}
Let us cap each of these tangles with a tangle $\langle \textrm{\textbf{G}}|$, which we write 
\be\label{G-state}
{\raisebox{-0.9cm}{\includegraphics[width=2cm]{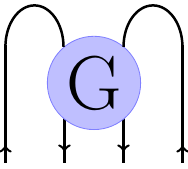}}} = \sum_{s,r_1,r_2} g_{s,r_1,r_2}\langle\phi^{(1)}_{s,r_1,r_2}(R,\bar{R},\bar{R},R)|~.
\ee
Then, the difference between the invariants of the mutant pairs arising from these 2-tangles will be
\begin{align}
{\raisebox{-1.8cm}{\includegraphics[width=5cm]{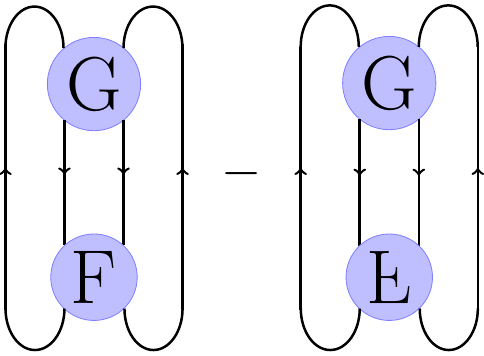}}}\ & = \ 2(f_{(1;1),0,1}~ g_{(1;1),0,1}+f_{(1;1),1,0}~ g_{(1;1),1,0})~, \\
{\raisebox{-1.8cm}{\includegraphics[width=5cm]{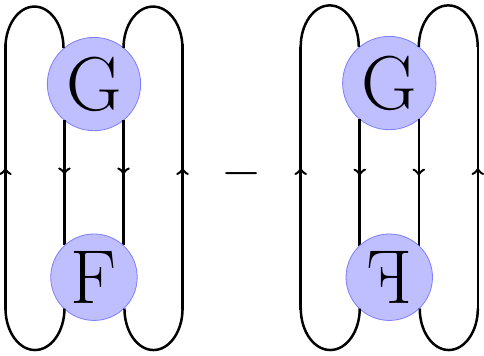}}} \ & = \ (f_{(1;1),0,1}+f_{(1;1),1,0})(g_{(1;1),0,1}+g_{(1;1),1,0})~.\label{diff-2-tangle}
\end{align}

\section{Example: Kinoshita-Terasaka and Conway}
In this section, we explicitly evaluate $\yng(2,1)$-colored HOMFLY-PT polynomials of Kinoshita-Terasaka $K_{KT}$ and Conway $K_C$ knot by means of braiding operations on conformal blocks. We refer the reader to  \cite{Gu:2014,Kaul:1991vt,Kaul:1992rs,Nawata:2013qpa,Ramadevi:1992dh} for technical details.

\begin{figure}[ht]
\center
\includegraphics[width=6cm]{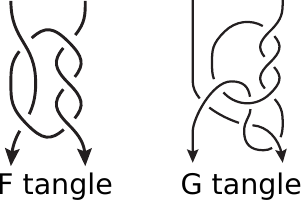}\caption{$\langle \textrm{\textbf{G}}|\textrm{\textbf{F}}\rangle$ and $\langle \textrm{\textbf{G}}|\protect\reflectbox{\textbf{F}}\rangle$ corresponds to Kinoshita-Terasaka $K_{KC}$ and Conway $K_C$ knot, respectively.}\label{F-G}
\end{figure}

\begin{figure}[ht]
\center
\subfloat[two boundaries]{\raisebox{.5cm}{\includegraphics[width=0.2\textwidth]{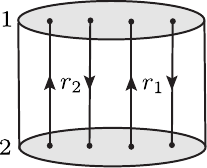}}} \qquad \qquad \qquad 
\subfloat[three boundaries]{\includegraphics[width=0.23\textwidth]{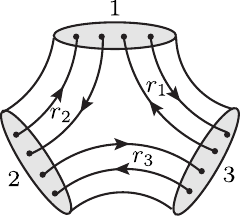}}
\caption{three-manifolds with boundaries }\label{fig:3bdry}
\end{figure}

For this mutant pair, the corresponding 2-tangles $\textrm{\textbf{F}}$ and $\textrm{\textbf{G}}$ are drawn in Figure \ref{F-G} and    the $M_y$-mutation is performed on the $\textrm{\textbf{F}}$-tangle. Although the difficulty lies on the explicit expressions for fusion matrices in this method, the fusion matrices for the representation $\yng(2,1)$ are delightedly  available thanks to the remarkable calculations in \cite{Gu:2014}. However, the quasi-plat representation Figure \ref{F-G} requires us to define the states incorporating multiplicity indices for the three-manifold with three $S^2$ boundaries shown in Figure \ref{fig:3bdry}(b).

First of all, $S^2\times I$ with four strands (Figure \ref{fig:3bdry}(a)) is equivalent to the identity operation. However, it can be also written by the tensor product of two ``ket'' states. In fact, bringing a ``bra'' to a ``ket'' amounts to the conjugation, which involves two $3j$-phases:
\bea
|\textrm{2-bdry}\rangle&=&\sum_{t,r_1,r_2}|\phi^{(1)}_{t;r_1,r_2}(R, \overline{R},R, \overline{R})\rangle~~
\langle\phi^{(1)}_{{t};r_1,r_2}(R, \overline{R},R, \overline{R})| \\
&=&\sum_{t,r_1,r_2}\{R,\bar{R},\bar{t},r_1\}\{R,\bar{R},\bar{t},r_2\}
|\phi^{(1)}_{t;r_1,r_2}(R, \overline{R},R, \overline{R})\rangle_{\textbf{1}}~~
|\phi^{(1)}_{{t};r_2,r_1}(R, \overline{R},R, \overline{R})\rangle_{\textbf{2}}~.\nonumber
\eea
Since capping the third boundary of Figure \ref{fig:3bdry}(b) leads to Figure \ref{fig:3bdry}(a), we have 
\bea
&&\sum_{t,r_1,r_3}\{R,\bar{R},\bar{t},r_3\}\delta_{r_1,r_3}\sqrt{\dim_q t}~~{}_{\textbf{3}}\langle\phi^{(1)}_{{t};r_3,r_1}(R, \overline{R},R, \overline{R}) |\textrm{3-bdry}\rangle \cr
&=&\sum_{t,r_1,r_2}\{R,\bar{R},\bar{t},r_1\}\{R,\bar{R},\bar{t},r_2\}
|\phi^{(1)}_{t;r_1,r_2}(R, \overline{R},R, \overline{R})\rangle_{\textbf{1}}~~
|\phi^{(1)}_{{t};r_2,r_1}(R, \overline{R},R, \overline{R})\rangle_{\textbf{2}}~.
\eea
Thus, the state $|\textrm{3-bdry}\rangle$ for Figure \ref{fig:3bdry}(b) can be defined as
\be\label{3-bdry}
|\textrm{3-bdry}\rangle=
\sum_{t,r_1,r_2,r_3}\Omega(t,r_1,r_2,r_3)
|\phi^{(1)}_{t;r_1,r_2}(R, \overline{R},R, \overline{R})\rangle_{\textbf{1}}
|\phi^{(1)}_{{t};r_2,r_3}(R, \overline{R},R, \overline{R})\rangle_{\textbf{2}}|\phi^{(1)}_{{t};r_3,r_1}(R, \overline{R},R, \overline{R})\rangle_{\textbf{3}}~,
\ee
where
$$
\Omega(t,r_1,r_2,r_3)=\frac{\{R,\bar{R},\bar{t},r_1\}\{R,\bar{R},\bar{t},r_2\}\{R,\bar{R},\bar{t},r_3\}}{\sqrt{\dim_{q}t}}~.
$$
Furthermore, it is straightforward to extrapolate it to multi-boundary states as 
$$
|n\textrm{-bdry}\rangle=\sum_{t,r_1,\cdots,r_n} \frac{\prod_{i=1}^n \{R,\bar{R},\bar{t},r_i\}}{\left(\sqrt{\dim_{q}t} \;\right)^{n-2}} \bigotimes_{i=1}^n |\phi^{(1)}_{t;r_i,r_{i+1}}(R, \overline{R},R, \overline{R})\rangle_{\textbf{i}}~.
$$

With these junctions with multi-boundaries, we can write down the states $|\textrm{\textbf{F}}\rangle$ in \eqref{F-state} and $|\textrm{\textbf{G}}\rangle$ in \eqref{G-state} by using the representations in Figure \ref{F-G}. With the help of the Racah backcoupling rule \eqref{identity}, the formulas can be simplified as follows:
\begin{align}\label{f-g-explicit}
f_{t,r_1,r_2}&=(\dim_q R)^2 \sum_{\cdots}\Omega(t,r_1,r_2,r_3)~(\lambda_{R,{R};wr_5}^{(+)})^{-3}a_{w;r_5,r_5}^{0;0,0}\begin{bmatrix}
 \bar{R}&R\\
R&\bar{R}
\end{bmatrix} a_{w;r_5,r_5}^{ t;r_2,r_3}\begin{bmatrix}
 \bar{R}&R\\
R& \bar{R}
\end{bmatrix}\cr
&\quad {(\lambda_{R,\bar{R};ur_4}^{(-)})^{2}a_{u;r_4,r_4}^{0;0,0}\begin{bmatrix}
 R&\bar{R}\\
R&\bar{R}
\end{bmatrix}a_{u;r_4,r_4}^{t;r_3,r_1}\begin{bmatrix}
 R&\bar{R}\\
R& \bar{R}
\end{bmatrix}^\ast}~,~ \cr
 g_{t,r_1, r_2}&=(\dim_q R)^2~\sum_{\cdots}~~
\Omega(i,\tilde{r}_1,\tilde{r}_2,\tilde{r}_3)\Omega(j,\tilde{r}_6,\tilde{r}_7,\tilde{r}_8)~(\lambda_{R,{R};l \tilde{r}_{5}}^{(+)})^{3}a_{l;\tilde{r}_5,\tilde{r}_5}^{0;0,0}\begin{bmatrix}
 \bar{R}&R\\
R&\bar{R}
\end{bmatrix} \nonumber\\
&\quad a_{l;\tilde{r}_5,\tilde{r}_5}^{ i;\tilde{r}_2,\tilde{r}_3}\begin{bmatrix}
 \bar{R}&R\\
R& \bar{R}
\end{bmatrix}^\ast
(\lambda_{R,{R};k \tilde{r}_{4}}^{(+)})^{-2}a_{k;\tilde{r}_4,\tilde{r}_4}^{0;0,0}\begin{bmatrix}
 \bar{R}&R\\
R&\bar{R}
\end{bmatrix}a_{k;\tilde{r}_4,\tilde{r}_4}^{ i;\tilde{r}_1,\tilde{r}_2}\begin{bmatrix}
 \bar{R}&R\\
R& \bar{R}
\end{bmatrix}
(\lambda_{R,\bar{R};s \tilde{r}_{9}}^{(-)})^{-2}\nonumber\\ 
&\quad a_{s;\tilde{r}_9,\tilde{r}_9}^{0;0,0}\begin{bmatrix}
 R&\bar{R}\\
R&\bar{R}
\end{bmatrix} 
 a_{s;\tilde{r}_9,\tilde{r}_9}^{ j;\tilde{r}_7,\tilde{r}_6}\begin{bmatrix}
 R&\bar{R}\\
R& \bar{R}
\end{bmatrix}^\ast (\lambda_{R,\bar{R};tr_1}^{(-)})a_{j;\tilde{r}_8,\tilde{r}_7}^{t; r_1,r_2}\begin{bmatrix}
R& \bar{R}\\
R&\bar{R}
\end{bmatrix}a_{j;\tilde{r}_8,\tilde{r}_6}^{i; \tilde{r}_{1},\tilde{r}_{3}}\begin{bmatrix}
R& \bar{R}\\
R&\bar{R}
\end{bmatrix}~.
\end{align}
Taking $t\equiv (1;1)$ and $ {r}_1 \neq  r_2 $, one can evaluate the difference \eqref{diff-2-tangle} of $\yng(2,1)$-colored HOMFLY-PT polynomials of the mutant pair by using the data of the quantum $6j$-symbols in \cite{Gu:2014}   
\begin{multline}\label{diff}
P_{\yng(2,1)}(K_{KT};a,q)-P_{\yng(2,1)}(K_{C};a,q)=a^{-5} q^{-18}(a-1)(a-q^2) (a q^2-1) (a-q^3)^2 (a
   q^3-1)^2\cr
(q-1)^2 (q^3-1)^2 (q^6-q^5+q^4-q^3+q^2-q+1)^2~.
\end{multline}
Then, it is easy to see from \eqref{diff} that the $\fraksl(2)$ ($a=q^2$) and $\fraksl(3)$ ($a=q^3$) quantum invariants as well as Alexander polynomials ($a=1$) cannot distinguish this mutant pair. 
The difference becomes apparent for $N>3$ and especially, at $N=4$, it factorizes as
$$
J^{\fraksl(4)}_{\yng(2,1)}(K_{KT};q)-J^{\fraksl(4)}_{\yng(2,1)}(K_{C};q)=-q^{-30}(1 - q)^6 (1 + q^2) (1 - q^3)^2 (1 - q^6) (1 - q^{14})^2~,
$$
which is consistent with the result obtained by Ochiai with the computer software ``Knot Theory By Computer'' \cite{Ochiai}.\footnote{It was programmed based on the cabling method \cite{Murakami:1989}. The result is present in \cite{Murakami:2000} and the change of variable $q\to q^2$ is necessary to see the equivalence.}

Furthermore, the explicit forms \eqref{f-g-explicit} allow us to evaluate colored HOMFLY-PT polynomials themselves. It turns out that this method is computationally efficient and it takes less than 15 minutes with a current desktop computer for the evaluation.\footnote{A mathematica file for this computation is linked on the arXiv page as an ancillary file.} The results are given in \eqref{HOMFLY-KT} and \eqref{HOMFLY-C}. Note that  the $\yng(2,1)$-colored HOMFLY-PT polynomial of a knot $K$ has the symmetry $P_{\yng(2,1)}(K;a,q)=P_{\yng(2,1)}(K;a,q^{-1})$ since the Young diagram $\yng(2,1)$ is invariant under transposition. Indeed, the invariants \eqref{HOMFLY-KT} and \eqref{HOMFLY-C} enjoy the symmetry.
At the $a=q^2$ specialization, they reduce to the Jones polynomial
$$
J^{\fraksl(2)}_{\yng(1)}(K_{KT};q)=J^{\fraksl(2)}_{\yng(1)}(K_{C};q)=\frac{1}{q^6}-\frac{2}{q^5}+\frac{2}{q^4}-\frac{2}{q^3}+\frac{1}{q^2}+2 q-2 q^2+2 q^3-q^4~.
$$

%
%
%

\newpage
\bea\label{HOMFLY-KT}
&&P_{\yng(2,1)}(K_{KT};a,q)\cr
&=&- a^3q^{-10} \Big(q^{20}-2 q^{19}+5 q^{18}-9 q^{17}+15 q^{16}-20 q^{15}+27 q^{14}-32
   q^{13}+38 q^{12}-40 q^{11}\cr
   &&\qquad\qquad+42 q^{10}-40 q^9+38 q^8-32 q^7+27 q^6-20 q^5+15 q^4-9
   q^3+5 q^2-2 q+1\Big)\cr
   &&+a^2q^{-13} \Big(q^{26}+2 q^{24}+q^{23}-4 q^{22}+19 q^{21}-32
   q^{20}+67 q^{19}-95 q^{18}+142 q^{17}-172 q^{16}\cr
   &&\qquad\qquad+218 q^{15}-228 q^{14}+246
   q^{13}-228 q^{12}+218 q^{11}-172 q^{10}+142 q^9-95 q^8\cr
   &&\qquad\qquad+67 q^7-32 q^6+19 q^5-4
   q^4+q^3+2 q^2+1\Big) \cr
   &&+aq^{-15} \Big(-2 q^{30}+2 q^{29}-9 q^{28}+12 q^{27}-27
   q^{26}+30 q^{25}-58 q^{24}+51 q^{23}-89 q^{22}\cr
   &&\qquad\qquad+85 q^{21}-133 q^{20}+102 q^{19}-163
   q^{18}+137 q^{17}-186 q^{16}+130 q^{15}-186 q^{14}\cr
   &&\qquad\qquad+137 q^{13}-163 q^{12}+102
   q^{11}-133 q^{10}+85 q^9-89 q^8+51 q^7-58 q^6+30 q^5\cr
   &&\qquad\qquad-27 q^4+12 q^3-9 q^2+2
   q-2\Big) \cr	
   &&+q^{-18}\Big(q^{35}+4 q^{33}-3 q^{32}+15 q^{31}-16 q^{30}+57 q^{29}-61
   q^{28}+131 q^{27}-142 q^{26}+248 q^{25}\cr
   &&\qquad\qquad-212 q^{24}+309 q^{23}-229 q^{22}+311
   q^{21}-170 q^{20}+263 q^{19}-141 q^{18}+263 q^{17}\cr
   &&\qquad\qquad-170 q^{16}+311 q^{15}-229
   q^{14}+309 q^{13}-212 q^{12}+248 q^{11}-142 q^{10}+131 q^9\cr
   &&\qquad\qquad-61 q^8+57 q^7-16 q^6+15
   q^5-3 q^4+4 q^3+q\Big) \cr
   &&+a^{-1}q^{-18}\Big(-q^{36}+q^{35}-6 q^{34}+8 q^{33}-24 q^{32}+22
   q^{31}-54 q^{30}+46 q^{29}-105 q^{28}+80 q^{27}\cr
   &&\qquad\qquad-185 q^{26}+168 q^{25}-347 q^{24}+332
   q^{23}-574 q^{22}+547 q^{21}-798 q^{20}+701 q^{19}\cr
   &&\qquad\qquad-888 q^{18}+701 q^{17}-798
   q^{16}+547 q^{15}-574 q^{14}+332 q^{13}-347 q^{12}+168 q^{11}\cr
   &&\qquad\qquad-185 q^{10}+80 q^9-105
   q^8+46 q^7-54 q^6+22 q^5-24 q^4+8 q^3-6 q^2+q-1\Big) \cr
   &&+a^{-2}q^{-18}\Big(q^{36}-q^{35}+6
   q^{34}-8 q^{33}+24 q^{32}-21 q^{31}+55 q^{30}-49 q^{29}+105 q^{28}-79 q^{27}\cr
   &&\qquad\qquad+183
   q^{26}-157 q^{25}+307 q^{24}-275 q^{23}+488 q^{22}-446 q^{21}+662 q^{20}-567
   q^{19}\cr
   &&\qquad\qquad+738 q^{18}-567 q^{17}+662 q^{16}-446 q^{15}+488 q^{14}-275 q^{13}+307
   q^{12}-157 q^{11}\cr
   &&\qquad\qquad+183 q^{10}-79 q^9+105 q^8-49 q^7+55 q^6-21 q^5+24 q^4-8 q^3+6
   q^2-q+1\Big) \cr
   &&-a^{-3}q^{-18}\Big(q^{35}+4 q^{33}-3 q^{32}+16 q^{31}-16 q^{30}+56 q^{29}-60
   q^{28}+130 q^{27}-144 q^{26}\cr
   &&\qquad\qquad+250 q^{25}-239 q^{24}+356 q^{23}-327 q^{22}+431
   q^{21}-351 q^{20}+452 q^{19}-368 q^{18}\cr
   &&\qquad\qquad+452 q^{17}-351 q^{16}+431 q^{15}-327
   q^{14}+356 q^{13}-239 q^{12}+250 q^{11}-144 q^{10}\cr
   &&\qquad\qquad+130 q^9-60 q^8+56 q^7-16 q^6+16
   q^5-3 q^4+4 q^3+q\Big) \cr
   &&+a^{-4}q^{-15} \Big(2 q^{30}-2 q^{29}+9 q^{28}-13 q^{27}+29
   q^{26}-30 q^{25}+57 q^{24}-52 q^{23}+84 q^{22}-77 q^{21}\cr
   &&\qquad\qquad+121 q^{20}-107 q^{19}+169
   q^{18}-165 q^{17}+213 q^{16}-176 q^{15}+213 q^{14}-165 q^{13}\cr
   &&\qquad\qquad+169 q^{12}-107
   q^{11}+121 q^{10}-77 q^9+84 q^8-52 q^7+57 q^6-30 q^5+29 q^4\cr
   &&\qquad\qquad-13 q^3+9 q^2-2
   q+2\Big) \cr
   &&-a^{-5}q^{-13}\Big(q^{26}+2 q^{24}-3 q^{22}+18 q^{21}-28 q^{20}+61
   q^{19}-94 q^{18}+144 q^{17}-178 q^{16}\cr
   &&\qquad\qquad+226 q^{15}-245 q^{14}+264 q^{13}-245
   q^{12}+226 q^{11}-178 q^{10}+144 q^9-94 q^8\cr
   &&\qquad\qquad+61 q^7-28 q^6+18 q^5-3 q^4+2
   q^2+1\Big) \cr
   &&+a^{-6}q^{-10} \Big(q^{20}-2 q^{19}+5 q^{18}-9 q^{17}+14 q^{16}-17 q^{15}+22
   q^{14}-25 q^{13}+29 q^{12}-29 q^{11}\cr
   &&\qquad\qquad+30 q^{10}-29 q^9+29 q^8-25 q^7+22 q^6-17 q^5+14
   q^4-9 q^3+5 q^2-2 q+1\Big) ~.\cr
   &&
   \eea
   \newpage
   \bea\label{HOMFLY-C}
&&P_{\yng(2,1)}(K_{C};a,q)\cr
&=&-a^3{q^{-10}}\Big(q^{20}-2 q^{19}+5 q^{18}-9 q^{17}+15 q^{16}-20 q^{15}+27 q^{14}-32 q^{13}+38
   q^{12}-40 q^{11}\cr
   &&\qquad\qquad+42 q^{10}-40 q^9+38 q^8-32 q^7+27 q^6-20 q^5+15 q^4-9 q^3+5 q^2-2
   q+1\Big)  \cr
   &&+a^2{q^{-13}}\Big(q^{26}+2 q^{24}+11 q^{21}-18 q^{20}+43 q^{19}-59
   q^{18}+93 q^{17}-110 q^{16}+146 q^{15}\cr
   &&\qquad\qquad-148 q^{14}+162 q^{13}-148 q^{12}+146 q^{11}-110
   q^{10}+93 q^9-59 q^8+43 q^7\cr
   &&\qquad\qquad-18 q^6+11 q^5+2 q^2+1\Big)  \cr
   &&+a{q^{-15}}\Big(-2
   q^{30}+2 q^{29}-7 q^{28}+5 q^{27}-15 q^{26}+11 q^{25}-28 q^{24}+12 q^{23}-43 q^{22}\cr
   &&\qquad\qquad+34
   q^{21}-85 q^{20}+59 q^{19}-125 q^{18}+110 q^{17}-166 q^{16}+110 q^{15}-166 q^{14}\cr
   &&\qquad\qquad+110
   q^{13}-125 q^{12}+59 q^{11}-85 q^{10}+34 q^9-43 q^8+12 q^7-28 q^6+11 q^5\cr
   &&\qquad\qquad-15 q^4+5
   q^3-7 q^2+2 q-2\Big)  \cr
   &&+{q^{-17}} \Big(q^{34}-q^{33}+6 q^{32}-3 q^{31}+11 q^{30}-5
   q^{29}+31 q^{28}-15 q^{27}+63 q^{26}-47 q^{25}\cr
   &&\qquad\qquad+130 q^{24}-81 q^{23}+169 q^{22}-91
   q^{21}+185 q^{20}-52 q^{19}+155 q^{18}-41 q^{17}\cr
   &&\qquad\qquad+155 q^{16}-52 q^{15}+185 q^{14}-91
   q^{13}+169 q^{12}-81 q^{11}+130 q^{10}-47 q^9\cr
   &&\qquad\qquad+63 q^8-15 q^7+31 q^6-5 q^5+11 q^4-3
   q^3+6 q^2-q+1\Big)\cr
   &&+{a^{-1}q^{-17} }\Big(-3 q^{34}+3 q^{33}-8 q^{32}+q^{31}-12 q^{30}-13 q^{29}-59
   q^{27}+38 q^{26}-144 q^{25}\cr
   &&\qquad\qquad+124 q^{24}-294 q^{23}+258 q^{22}-471 q^{21}+411 q^{20}-628
   q^{19}+509 q^{18}-690 q^{17}\cr
   &&\qquad\qquad+509 q^{16}-628 q^{15}+411 q^{14}-471 q^{13}+258
   q^{12}-294 q^{11}+124 q^{10}-144 q^9\cr
   &&\qquad\qquad+38 q^8-59 q^7-13 q^5-12 q^4+q^3-8 q^2+3
   q-3\Big) \cr
   &&+{a^{-2}q^{-17} }\Big(3 q^{34}-3 q^{33}+8 q^{32}-q^{31}+13 q^{30}+14 q^{29}-3 q^{28}+59
   q^{27}-37 q^{26}+142 q^{25}\cr
   &&\qquad\qquad-113 q^{24}+254 q^{23}-201 q^{22}+385 q^{21}-310 q^{20}+492
   q^{19}-375 q^{18}+540 q^{17}\cr
   &&\qquad\qquad-375 q^{16}+492 q^{15}-310 q^{14}+385 q^{13}-201
   q^{12}+254 q^{11}-113 q^{10}+142 q^9\cr
   &&\qquad\qquad-37 q^8+59 q^7-3 q^6+14 q^5+13 q^4-q^3+8 q^2-3
   q+3\Big) \cr
   &&+{a^{-3}q^{-17} }\Big(-q^{34}+q^{33}-6 q^{32}+3 q^{31}-12 q^{30}+5 q^{29}-30
   q^{28}+14 q^{27}-62 q^{26}\cr
   &&\qquad\qquad+49 q^{25}-132 q^{24}+108 q^{23}-216 q^{22}+189 q^{21}-305
   q^{20}+233 q^{19}-344 q^{18}\cr
   &&\qquad\qquad+268 q^{17}-344 q^{16}+233 q^{15}-305 q^{14}+189
   q^{13}-216 q^{12}+108 q^{11}-132 q^{10}\cr
   &&\qquad\qquad+49 q^9-62 q^8+14 q^7-30 q^6+5 q^5-12 q^4+3
   q^3-6 q^2+q-1\Big) \cr
   &&+{a^{-4}q^{-15}}\Big(2 q^{30}-2 q^{29}+7 q^{28}-6 q^{27}+17 q^{26}-11
   q^{25}+27 q^{24}-13 q^{23}+38 q^{22}\cr
   &&\qquad\qquad-26 q^{21}+73 q^{20}-64 q^{19}+131 q^{18}-138
   q^{17}+193 q^{16}-156 q^{15}+193 q^{14}\cr
   &&\qquad\qquad-138 q^{13}+131 q^{12}-64 q^{11}+73 q^{10}-26
   q^9+38 q^8-13 q^7+27 q^6-11 q^5\cr
   &&\qquad\qquad+17 q^4-6 q^3+7 q^2-2 q+2\Big) \cr
   &&-{a^{-5}q^{-13} } \left(q^2+q+1\right)^2 \Big(q^{22}-2 q^{21}+3 q^{20}-3 q^{19}+q^{18}+13
   q^{17}-40 q^{16}+79 q^{15}\cr
   &&\qquad\qquad-123 q^{14}+171 q^{13}-207 q^{12}+222 q^{11}-207 q^{10}+171
   q^9-123 q^8+79 q^7\cr
   &&\qquad\qquad-40 q^6+13 q^5+q^4-3 q^3+3 q^2-2 q+1\Big)\cr
   &&+a^{-6}q^{-10} \Big(q^{20}-2 q^{19}+5 q^{18}-9 q^{17}+14 q^{16}-17 q^{15}+22 q^{14}-25
   q^{13}+29 q^{12}-29 q^{11}\cr
   &&\qquad\qquad+30 q^{10}-29 q^9+29 q^8-25 q^7+22 q^6-17 q^5+14 q^4-9 q^3+5
   q^2-2 q+1\Big)~.\cr
   &&
\eea

\newpage

As a concluding remark, we should mention that it is easy to find infinitely many mutant pairs that can be distinguished by $\yng(2,1)$-colored HOMFLY-PT polynomials. For instance, we can add even number of crossings to some part of Figure \ref{F-G}. Furthermore, as we have illustrated in section 2,  $\yng(2,1)$-colored HOMFLY-PT polynomials can detect a mutation on a two-tangle $|\textrm{\textbf{F}}\rangle$ with $\langle\phi^{(1)}_{t,r_1,r_2}(R,\bar R,R,\bar R)|\textrm{\textbf{F}}\rangle\neq0$ for $(t,r_1,r_2)=((1;1),0,1)$ or $((1;1),1,0)$, and asymmetric two-tangles generically obey this condition. The forthcoming paper \cite{Mironov:2015zz} will give examples of mutant pairs which are distinguishable by $\yng(2,1)$-colored HOMFLY-PT polynomials and those which are indistinguishable.  In addition,  it would be interesting to investigate the conjecture \cite{Gu:2014yba} that any mutant pair can be distinguished by some HOMFLY-PT polynomials colored by Young diagrams with two rows.

As we have seen, properties of colored HOMFY polynomials with multiplicity structure are very different from those in the multiplicity-free cases. The study of HOMFLY-PT polynomials colored by non-rectangular Young tableaux are still very limited \cite{Anokhina:2012rm,Anokhina:2013ica,Gu:2014,Mironov:2014zza} and the properties have not been fully uncovered yet. Presumably, intricate feature behind multiplicity will  be reflected to higher rank Witten-Reshetikhin-Turaev invariants of three-manifolds, higher rank volume conjectures, knot homology colored by general Young diagrams.  We hope
that our results will serve as a stepping stone towards the study of higher rank quantum and homological invariants of knots and three-manifolds.

\section*{Acknowledgement}
The authors would like to express sincere gratitude to H. Morton for correspondence. In addition, they are grateful to J. Gu, H. Jockers, R.K. Kaul, A. Morozov, A. Mironov and Zodinmawia for discussion and correspondence. They are also obliged to S. Garoufalidis, Hitoshi Murakami, Jun Murakami and P. Su{\l}kowski for comments on the manuscript. The comments from the referee of Journal of Knot Theory and Its Ramifications improves the content so that we would like to thank the referee. S.N. would like to thank American Institute of Mathematics, conference and workshop ``Group Theory and Knots'' at Natal, University of Bonn for the warm hospitality where a part of the work had been carried out. The work of S.N. is supported by the ERC Starting Grant no.~335739 \textit{``Quantum fields and knot homologies''}, funded by the European Research Council under the European Union's Seventh Framework Programme.

\bibliographystyle{alpha}

\end{document}